\documentclass[12pt,a4article]{article}

\usepackage{tikz}
\usepackage{geometry}
\usepackage{amsmath,amsfonts}
\usepackage{amssymb}
\usepackage{array}
\usepackage{mathrsfs}
\usepackage{latexsym}
\usepackage{xcolor,graphicx,float,graphics}
\usepackage{subfigure}
\usepackage{subcaption}
\usepackage{bm}
\usepackage{multirow}
\usepackage{epstopdf}
\usepackage{pst-poly}     
\usepackage{cite}
\usepackage{caption}
\usepackage{setspace}
\usepackage[modulo, pagewise, left, mathlines]{lineno}
\newtheorem{thm}{Theorem}[section]
\newtheorem{lem}[thm]{Lemma}
\newtheorem{Def}{Definition}[section]
\newtheorem{cor}{Corollary}[section]
\newtheorem{prop}{Proposition}[section]

\newtheorem{rem}{Remark}[section]

\newtheorem{cla}[thm]{Claim}

\newtheorem{fa}{Fact}[section]
\makeatletter \@addtoreset{equation}{section} \makeatother

\begin{document}

\title{Fan's condition for completely independent spanning trees}
\author{Jie Ma$^a$, Junqing Cai$^{a,b}$\thanks{Corresponding author: Junqing Cai (caijq09@163.com)}\\
\mbox{}\hspace{0.15cm}{\scriptsize $^a$ School of Mathematical Science, Tianjin Normal University, Tianjin, China}\\
{\scriptsize $^b$ Institute of Mathematics and Interdisciplinary Sciences, Tianjin Normal University, Tianjin, China}\\
\mbox{}\hspace{0.15cm}{\scriptsize E-mails: majielx@163.com, caijq09@163.com}}
\date{}
	
\maketitle

\begin{abstract}
Spanning trees $T_1,T_2, \dots,T_k$ of $G$ are $k$ completely independent spanning trees if, for any two vertices $u,v\in V(G)$, the paths from $u$ to $v$ in these $k$ trees are pairwise edge-disjoint and internal vertex-disjoint. Hasunuma proved that determining whether a graph contains $k$ completely independent spanning trees is NP-complete, even for $k = 2$. Araki posed the question of whether certain known sufficient conditions for hamiltonian cycles are also also guarantee two completely independent spanning trees? In this paper, we affirmatively answer this question for the Fan-type condition. Precisely, we proved that if $G$ is a connected graph such that each pair of vertices at distance 2 has degree sum at least $|V(G)|$, then $G$ has two completely independent spanning trees.
\vspace{1mm}

\end{abstract}

\section{Introduction}
	
In this paper, we consider only undirected, finite and simple graphs. For a graph $G$, we denote its vertex set and edge set by $V(G)$ and $E(G)$, respectively. The number of vertices in $G$ is called the order of $G$, written as $|V(G)|$ or $|G|$. The neighborhood and degree of a vertex $v$ in $G$ are the set and number of vertices which are adjacent to $v$ in $G$, denoted by $N_G(v)$ and $d_G(v)$, respectively. The minimum degree of $G$ is  $\delta (G)$=min$\{d(v)|v \in V(G)\}$. For a subgraph $H$ of $G$ and a subset $U$ of $V(G)$, let $N_H(v)=N_G(v)\cap V(H)$, $d_H(v)=|N_H(v)|$ and $N_H(U)=\cup_{u\in U} N_H(u)$. When there is no confusion, we simplify the notation to $N(v)$, $d(v)$ and $N(U)$. If $G$ is not a complete graph, then we define
$$\sigma_2(G)=\min\{d(u)+d(v)|u,v \in V(G)~and~uv\notin E(G)\},$$
$$\mu_2(G)=\min\{d(u)+d(v)|u,v \in V(G)~ and ~d(u,v)=2\}.$$
Otherwise, set $\sigma_2(G)=\mu_2(G)=+\infty$. 
The subgraph induced by $U$ is denoted by $G[U]$, and the subgraph induced by $V(G) \setminus U$ is denoted by $G-U$.  For a single vertex $u$, we use $G-u$ instead of $G-\{u\}$. The set of integers from 1 to $n$  is denoted by $[n]$.

A tree $T$ of a graph $G$ is a spanning tree if $V(T)$ = $V(G)$. In 2001, Hasunuma \cite{cist's concept} introduced the concept of completely independent spanning trees (CISTs).

\begin{Def}\emph{(Hasunuma \cite{cist's concept})}
Let $T_1,T_2, \dots,T_k$ be spanning trees in a graph $G$. If for any two vertices $u,v$ of $G$, the paths from $u$ to $v$ in $T_1,T_2, \dots,T_k$ are pairwise edge-disjoint and internal vertex-disjoint, then $T_1,T_2, \dots,T_k$ are $k$ completely independent spanning trees $(k$ CISTs).
\end{Def}

Completely independent spanning trees are a powerful tool for enhancing the reliability, redundancy, and efficiency of various systems, particularly those that rely on network structures. Hasunuma \cite{cist-NP} proved that determining whether a given graph contains $k$ CISTs is NP-complete, even for $k=2$. 

Let $V_1, V_2, \dots, V_k$ be $k$ disjoint sets of vertices in $G$. If $V(G)=V_1\cup V_2\cup \dots \cup V_k$, then we call $(V_1, V_2, \dots, V_k)$ a $k$-partition of $V(G)$. For $i\neq j$, $B(V_i, V_j, G)$ denotes the bipartite subgraph of $G$ induced by the edges of $G$ with one endvertex in $V_i$ and the other in $V_j$. If there is no confusion, we use $B(V_i, V_j)$ instead of $B(V_i, V_j, G)$.

\begin{Def}\emph{(Araki \cite{cist partition})}
Let $G$ be a graph and $(V_1, V_2, \dots, V_k)$ be a $k$-partition of $V(G)$. Then $(V_1,V_2, \dots, V_k)$ is a $k$-CIST-partition if it satisfies the following two conditions:
	
\emph{(i)} for each $i\in [k]$, the induced subgraph $G[V_i]$ is connected; and
	
\emph{(ii)} for any $1\le i\neq j\le k$, $B(V_i, V_j)$ has no tree component (i.e. each component of $B(V_i, V_j)$ has a cycle).
\end{Def}

Araki \cite{cist partition} gave an useful equivalent characterization of $k$ CISTs, which is crucial for proving their existence.

\begin{thm} \emph{(Araki \cite{cist partition}) }\label{thm1}
A connected graph has $k$ CISTs if and only if it has a $k$-CIST-partition.
\end{thm}

In the same paper, Araki \cite{cist partition} also showed that two well-established conditions for hamiltonian cycles-Dirac's and a condition obtained by Fleischner \cite{Fleischner hamlition}, are sufficient for the existence of two CISTs. 

\begin{thm}\emph{(Araki \cite{cist partition}) }
Let $G$ be a graph of order $n \ge 7$. If $\delta(G) \ge \frac{n}{2} $, then $G$ has two CISTs.	
\end{thm}

\begin{thm}\emph{(Araki \cite{cist partition}) }
	Let $G$ be a $2$-connected graph of order $n \ge 4$. Then the square $G^2$ has two CISTs.
\end{thm}

Based on these results, Araki \cite{cist partition} posed the question of whether certain known sufficient conditions for hamiltonian cycles are also also guarantee two completely independent spanning trees? Fan, Hong and Liu \cite{fan cist} confirmed this for Ore's condition.

\begin{thm} \emph{(Fan, Hong and Liu\cite{fan cist}) }
Let $G$ be a graph of order $n \ge 8$. If $\sigma_2(G)\ge n$, then $G$ has two CISTs.
\end{thm}

Hong and Zhang \cite{hong xia} extended this result to a condition by Gould and Jacobson \cite{hong-ha}.

\begin{thm}\emph{(Hong and Zhang \cite{hong xia})}\label{th3}
Let $G$ be a graph of order $n\ge 5$. If $|N(x) \cup N(y)| \ge \frac{n}{2}$ and $|N(x) \cap N(y)| \ge 3$ for every pair of nonadjacent vertices $x$ and $y$, then $G$ has two CISTs. (Note: Qin et al.\cite{hao} later pointed that $n \ge 8$ is the correct lower bound.)
\end{thm}

Fan's condition is another well-known condition for hamiltonian cycles.

\begin{thm}\emph{(Fan \cite{fan hamiltion})}\label{th6}
Let $G$ be a $2$-connected graph of order $n \ge 3$. If \textnormal{max}$\{d(u),d(v)\} \ge \frac{n}{2}$ for every pair of vertices $u$ and $v$ with $d(u,v)=2$, then $G$ has a hamiltonian cycle.
\end{thm}

From Theorem \ref{th6}, we derive the following corollary. 

\begin{cor}
Let $G$ be a connected graph of order $n \ge 3$. If $\mu_2(G) \ge n$, then $G$ has a hamiltonian cycle.
\end{cor}

In this paper, we prove that the condition $\mu_2(G) \ge n$ is also sufficient for the existence of two CISTs.

\begin{thm}\label{th2}
	Let $G$ be a connected graph of order $n\geq 7$. If $\mu_2(G) \ge n$, then $G$ has two CISTs.
\end{thm}

\begin{rem}	
$(1)$ The complete bipartite graph $K_{3,3}$ implies that $n\ge 7$ in Theorem \ref{th2} is necessary.\\
$(2)$ Let $G_1=K_s$ and $G_2=K_t$ be two disjoint complete graphs. The graph $G$ obtained from two disjoint complete graphs $K_s$ and $K_t$ by adding a new vertex $u$ and edges $\{uv|v\in V(K_s\cup K_t)\}$ implies that the condition $\mu_2(G) \ge n$ of Theorem \ref{th2} is sharp. 
\end{rem}

\section{Preliminary}

The following lemmas are useful in the proof of Theorem \ref{th2}.

\begin{lem}\label{lem1}
Let $G$ be a connected non-complete graph of order $n\geq 3$ and $\mu_2(G)\geq n$. Then for any two vertices $x$, $y$ with $d(x,y)=2$, $|N(x) \cap N(y)| \ge 2$. Furthermore, if $|N(x) \cap N(y)|=2$, then $ N(x) \cup N(y) = V(G) \setminus \{x, y\}$.
\end{lem}	

\noindent
\textbf{Proof.}  Note that $N(x) \cup N(y) \subseteq V(G) \setminus \{x,y\}$, we have $n\leq d(x)+d(y)=|N(x) \cup N(y)|+|N(x) \cap N(y)|\le |V(G) \setminus \{x,y\}|+|N(x) \cap N(y)|=(n-2)+|N(x) \cap N(y)|.$ Thus $|N(x) \cap N(y)| \ge 2$, and $|N(x) \cap N(y)|=2\Rightarrow N(x) \cup N(y) = V(G) \setminus \{x, y\}$. $\hfill\Box$

\begin{lem}\label{lem2}
Let $G$ be a connected graph of order $n\ge 3$. If $\mu_2(G)\geq n$, then
	
\textnormal{(i)} $G$ is $2$-connected$;$ 
	
\textnormal{(ii)} if $\{u,v\}$ is a cut of $G$, then $G-\{u,v\}$ has exactly two components.
\end{lem}

\noindent
\textbf{Proof.} (i) Assume for contradiction that $G$ is not 2-connected. Then there is a cut-vertex $u$ in $G$. Let $G_1$ and $G_2$ be two components of $G-u$. Since $G$ is connected, there exist vertices $v_1 \in V(G_1)$ and $v_2 \in V(G_2)$ such that $d(v_1,v_2)=2$. Thus, $n \le d(v_1)+d(v_2) \le |V(G_1)|+|V(G_2)|<n$, a contradiction. Hence, $G$ is $2$-connected.

(ii) Assume for contradiction that $\{u,v\}$ is a vertex cut of $G$ such that $G-\{u,v\}$ has three components $G_1$, $G_2$ and $G_3$. Then there exist vertices $v_1 \in V(G_1)$ and $v_2 \in V(G_2)$ such that $d(v_1,v_2)=2$. Thus, $n \le d(v_1)+d(v_2) \le (|V(G_1)|+1)+(|V(G_2)|+1)<n$, a contradiction. So $G-\{u,v\}$ has exactly two components. $\hfill\Box$

\begin{lem}\label{lem3}
Let $G$ be a connected graph of order $n\ge 7$ and $\mu_2(G)\geq n$. If the connectivity $\kappa(G)=2$, then $G$ has two CISTs.
\end{lem}

\noindent
\textbf{Proof.} Let $\{u,v\}$ be a vertex cut of $G$. Then by Lemma \ref{lem2}, $G-\{u,v\}$ has exactly two components $G_1$ and $G_2$. Without loss of generality, assume $|V(G_1)| \le |V(G_2)|$. 

Since $G$ is connected, there are vertices $x_1, x_2\in V(G_1)$  and $y_1, y_2\in V(G_2)$ (maybe $x_1=x_2$) such that $d(x_1,y_1)=d(x_2,y_2)=2$. By Lemma \ref{lem1}, $N(x_i) \cap N(y_i) = \{u,v\}$, $N(x_i)=(V(G_1) \setminus \{x_i\}) \cup \{u,v\}$ and $N(y_i)=(V(G_2) \setminus \{y_i\}) \cup \{u,v\}$ for  $i\in [2]$. 

If $|V(G_1)| \ge 2$, then let $V_1=\{u,y_1\} \cup V(G_1) \setminus\{x_1\}$, $V_2=\{v,x_1\} \cup V(G_2) \setminus\{y_1\}$. Obviously, $G[V_1]$, $G[V_2]$  and $B(V_1,V_2)$ are connected, and $ux_1x_2vy_1y_2u$ is a cycle of $B(V_1,V_2)$ (see Figure 1 (a)). Therefore, $(V_1,V_2)$ is 2-CIST-partition of $G$. 
\begin{figure}[htbp]
    \centering
    \begin{minipage}[b]{0.3\textwidth}
        \centering
        \includegraphics[width=\textwidth]{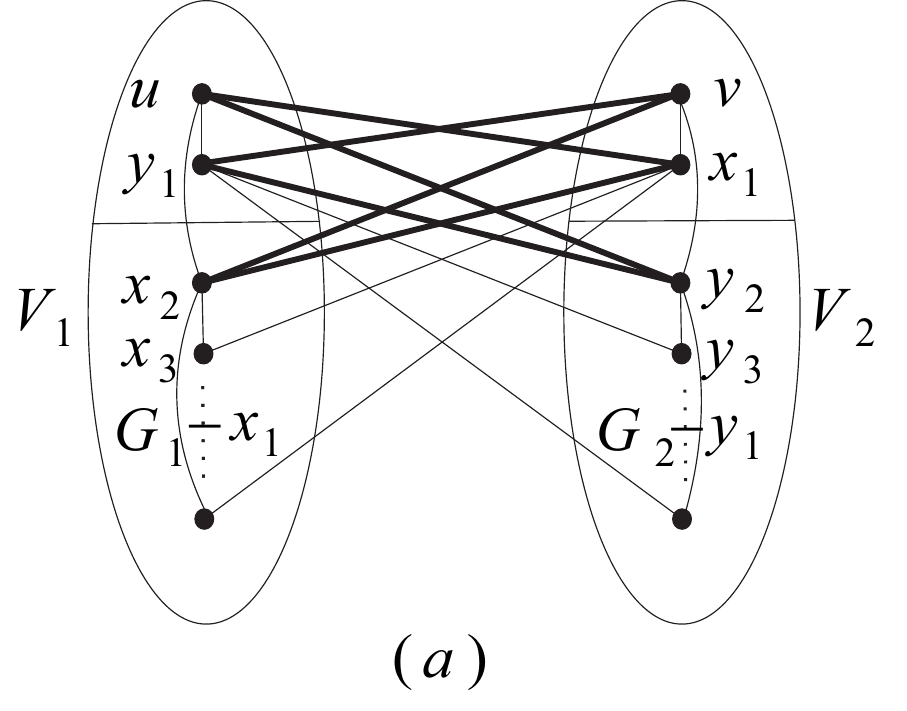}
        \label{fig:sub1}
    \end{minipage}
    \hfill
    \begin{minipage}[b]{0.3\textwidth}
        \centering
        \includegraphics[width=\textwidth]{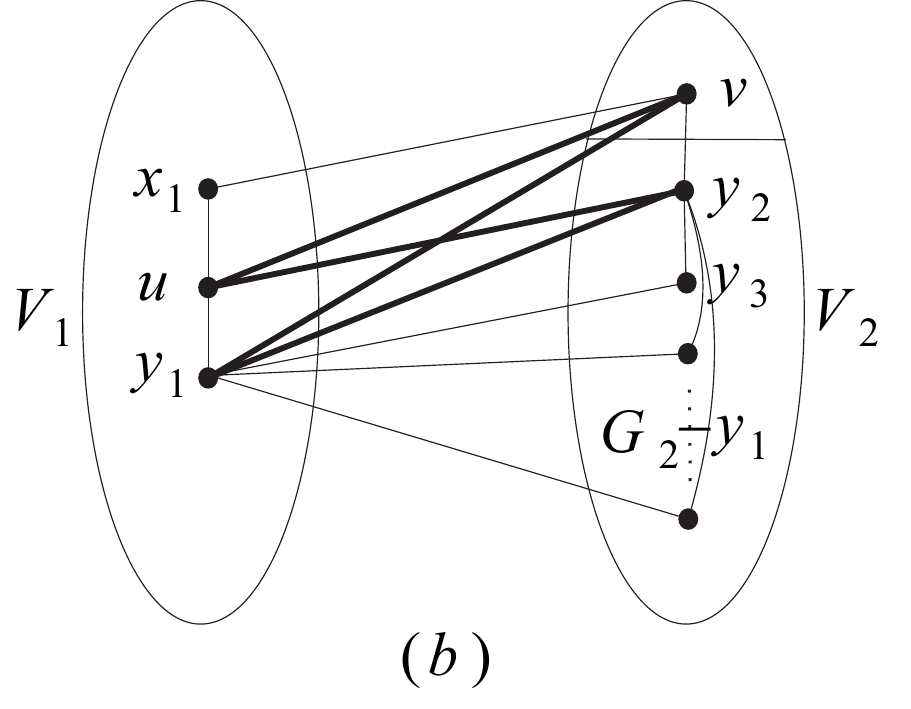}
        \label{fig:sub2}
    \end{minipage}
    \hfill
    \begin{minipage}[b]{0.3\textwidth}
        \centering
        \includegraphics[width=\textwidth]{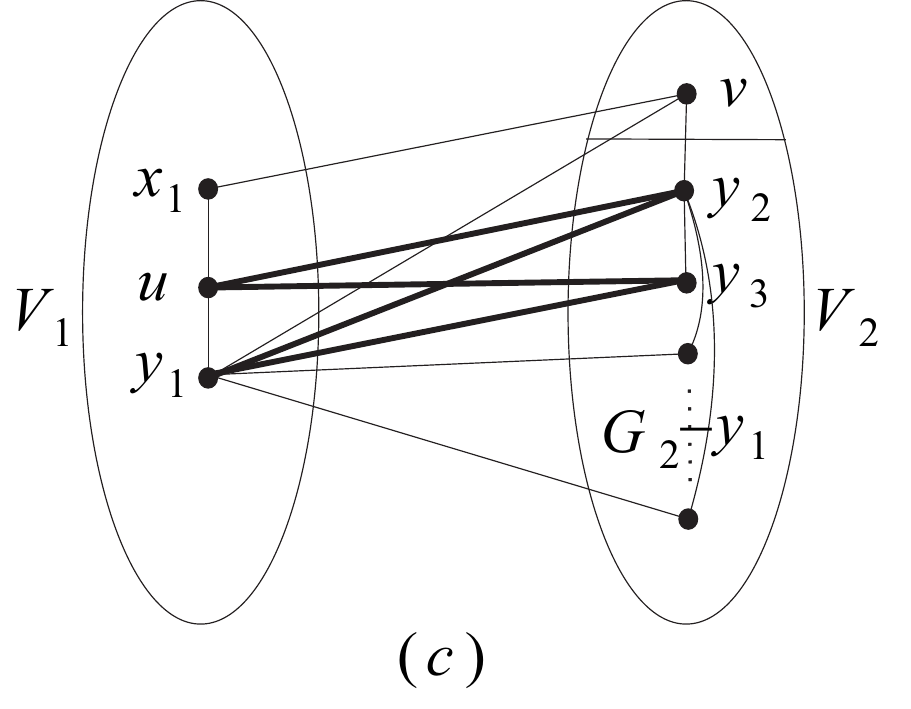}
        \label{fig:sub3}
    \end{minipage}
    \vspace{-0.5cm}
    \caption{\footnotesize The illustration of a 2-CIST-partition in Lemma \ref{lem3}.}
    \label{fig:three_figures}
\end{figure}

If $|V(G_1)| =1$, then $|V(G_2)|\geq 4$. Let $V_1=\{u,x_1,y_1\}$ and $V_2=\{v\} \cup V(G_2) \setminus\{y_1\}$. Obviously, $G[V_1]$, $G[V_2]$ and $B(V_1,V_2)$ are connected. If $uv \in E(G)$, then $uy_2y_1vu$ is a cycle of $B(V_1,V_2)$ (see Figure 1 (b)). If $uv \notin  E(G)$, then $d(u,v)=2$ and $d(u)+d(v) \ge n \ge 7$. So, there must exist a vertex $y_3 \in V(G_2)$ such that $uy_3 \in E(G)$ or $vy_3 \in E(G)$. In both cases, $d(x_1,y_3)=2$. Then by Lemma \ref{lem1}, $N(y_3)=(V(G_2) \setminus \{y_3\}) \cup \{u,v\}$. Thus $uy_2y_1y_3u$ is a cycle of $B(V_1,V_2)$ (see Figure 1 (c)). Therefore, $(V_1,V_2)$ is a 2-CIST-partition of $G$. 
 
Therefore, by Theorem \ref{thm1}, $G$ has two CISTs. 
$\hfill\Box$
	
\section{Proof of Theorem \ref{th2}}

Let $G$ be a graph satisfying the condition of Theorem \ref{th2}. By Lemma \ref{lem2}, $\kappa(G)\geq 2$. If $G$ is a complete graph, then $G$ has two CISTs. If $\kappa(G)=2$, then by Lemma \ref{lem3}, $G$ has two CISTs.

Next, we assume that $G$ is not a complete graph and $\kappa(G)\geq 3$. Then there are two vertices $x,y\in V(G)$ such that $d(x,y)= 2$. We choose such vertices $x$ and $y$ satisfying the following conditions:

(C1) $|N(x) \cap N(y)|$ is as small as possible;

(C2) $d(x)+d(y)$ is  as large as possible, subject to (C1).

Without loss of generality, suppose $d(x) \le d(y)$. Partition $V(G)\setminus \{x,y\}$ as follows:
\begin{align*}  
	M &= N(x) \cap N(y),\\
	X &= N(x)\setminus M,\\
	Y &= N(y)\setminus M,\\
	D &= V(G)\setminus (\{x,y\}\cup M \cup X\cup Y ).
\end{align*}
Then by Lemma \ref{lem1}, we have $|X|\leq |Y|$ and $|M| \ge 2$ . Since $n \le d(x)+d(y) \le |V(G)\setminus \{x,y\}|+|M|-|D|=n-2+|M|-|D|$, we have $|D|\leq |M|-2$. For convenience, set
$$X=\{x_1,x_2,\dots ,x_t\},~Y=\{y_1,y_2,\dots ,y_s\},~M=\{u_1,u_2,\dots ,u_m\}~\text{and}~D=\{v_1,v_2,\dots ,v_d\}.$$ 

By the choice of $x$ and $y$, for each pair of vertices $u,v\in V(G)$ with $d(u,v)=2$, we have $|N(u)\cap N(v)|\geq m$ and $n\leq d(u)+d(v)\leq d(x)+d(y)$.

\begin{cla}\label{cla7}
If $D\neq \varnothing$, then each vertex of $D$ has at least one neighbor in $M\cup Y$.
\end{cla}

\noindent
\textbf{Proof.} Suppose not. Let $v_i\in D$ be a vertex such that $N_{ M\cup Y} (v_i)=\varnothing $. Let $D'$ be a component of $G[D]$ containing $v_i$. We can claim that for each vertex $v_j\in D$, if $N_{ M\cup Y} (v_j)=\varnothing $, then $N_{X} (v_j)=\varnothing $, which implies that $N(v_j)\subseteq D$. Otherwise, $d(x,v_j)=2$. Hence, $n \le d(x)+d(v_j) \le (|X|+|M|)+(|X|+D|-1)<n$, a contradiction.  

Since $\kappa(G)\geq 3$, there is a vertex $v_j \in V(D')$ such that $N_{M\cup Y} (v_j)\neq \varnothing$. Let $P$ be a path in $D'$ connecting $v_i$ and $v_j$. Then there are two vertices $v_l, v_k\in V(P)$ such that $v_lv_k\in E(P)$, $N_{M\cup Y} (v_l)\neq \varnothing $ and $N_{M\cup Y} (v_k)= \varnothing $. Let $z\in N_{M\cup Y} (v_l)$. Then $d(z,v_k)=2$ and $|N(z)\cap N(v_k)| \ge m$. However, $N(v_k) \subseteq D$, this contradicts to  $m\ge |D|+2$. $\hfill\Box$

\begin{cla}\label{3-2 com}
If $\kappa(G)=3$, then for every cut $\{u,v,w\}$ of $G$,  $G-\{u,v,w\}$ has exactly two components.
\end{cla}

\noindent
\textbf{Proof.} Suppose not. Let $\{u,v,w\}$ be a cut of $G$ such that $G-\{u,v,w\}$ has at least three components. Suppose that $G_1$, $G_2$ and $G_3$ are three components of $G-\{u,v,w\}$ with $|V(G_1)| \le |V(G_2)| \le |V(G_3)|$. Since $\{u,v,w\}$ is a minimal cut, there are two vertices $z_1 \in V(G_1)$ and $z_2 \in V(G_2)$ such that $d(z_1,z_2) = 2$. Then, $n \le d(z_1)+d(z_2) \le |V(G_1)|+|V(G_2)|+4=n+1-|V(G_3)|$. So $|V(G_3)| \le 1$. So $|V(G_1)|=|V(G_2)|=|V(G_3)|=1$. Therefore, $7\leq n\le d(z_1)+d(z_2) \le 6$, a contradiction. $\hfill\Box$

\vspace{0.50cm}
\noindent{\bf Case 1}  $|M|=2$.

\vspace{0.5cm}
In this case, we have $M=\{u_1,u_2\}$ and $D = \varnothing$. Since $d(x) \ge \delta(G) \ge \kappa(G) \ge 3$, $X \ne \varnothing$.
	
\begin{cla}\label{de4}
$d(u_1) \ge 4$ or $d(u_2) \ge 4$.
\end{cla}	
	
\noindent
\textbf{Proof.} Suppose not. Then $d(u_1)=d(u_2) = 3$. If $u_1u_2 \in E(G)$, then $N(u_1)=\{x,y,u_2\}$ and $N(u_2)=\{x,y,u_1\}$. Then $d(x_1,u_1)=2$. However, $N(u_1) \cap N(x_1)=\{x\}$, this contradicts to Lemma \ref{lem1}. If $u_1u_2 \notin E(G)$, then $d(u_1,u_2)=2$. Thus, $6=d(u_1) +d(u_2) \ge n \ge 7$, a contradiction. $\hfill\Box$

\vspace{0.50cm}
By Claim \ref{de4}, we can assume without loss of generality that $d(u_1) \ge 4$.
	
\begin{cla}\label{cla3.3}
If  $d(u_2) \ge 4$ or $u_1u_2 \notin E(G)$, then $\{x,y,u_1\}$ is not a cut of $G$.
\end{cla}
	
\noindent
\textbf{Proof.} Suppose not. Then, by Claim \ref{3-2 com}, $G-\{x,y,u_1\}$ has exactly two components $G_1$ and $G_2$. We can claim that all the vertices of $X$ are in the same component. Otherwise, there is a vertex $x_i \in V(G_i)$ for $i\in[2]$ such that $d(x_1,x_2) = 2$. By Lemma \ref{lem1}, $N(x_1) \cap N(x_2) = \{x,u_1\}$, $N(x_1) \cup N(x_2) = V(G) \setminus \{x_1,x_2\}$. However, $y \notin N(x_1) \cup N(x_2)$, a contradiction. Thus, all the vertices of $X$ are in the same component of $G-\{x,y,u_1\}$. Similarly, all the vertices of $Y$ are in the same component of $G-\{x,y,u_1\}$. 

Next, we will prove that $X$ and $Y$ are contained in the same component of $G-\{x,y,u_1\}$. Otherwise, without loss of generality, assume $X\subseteq V(G_1)$ and $Y\subseteq V(G_2)$. Then $G -\{u_1,y\}$ is disconnected if $u_2\in V(G_1)$ or $G -\{u_1,x\}$ is disconnected if $u_2\in V(G_2)$. This contradicts to the assumption $\kappa(G) \ge 3$. Therefore, $X \cup Y$ is in the same component. Suppose without loss of generality that $V(G_1)=X \cup Y$ and $V(G_2)=\{u_2\}$. Thus $d(u_2) = 3$ and $u_1u_2\in E(G)$, which contradicts to the condition of Claim \ref{cla3.3}. Therefore, $\{x,y,u_1\}$ is not a cut of $G$. $\hfill\Box$

\vspace{0.5cm}
\noindent{\bf Subcase 1.1}  $d(u_2) \ge 4$.

Let $V_{1}=\{x,y,u_1\}$ and $V_{2}=X \cup Y \cup \{u_2\}$.
Obviously, $G[V_{1}]$ is connected. By Claim \ref{cla3.3},  $G[V_{2}]$ is connected. Since $d(u_1) \ge 4$, $B(V_{1},V_{2})$ is connected and $|N(u_1)\cap V_2|\geq 2$.

If $|N(u_1) \cap X| \ge 2$, without loss of generality, we assume $\{x_1,x_2\} \subseteq N(u_1)$, then $xx_1u_1x_2x$ is a cycle of $B(V_{1},V_{2})$. Similarly, if $|N(u_1) \cap Y| \ge 2$, then we can get a cycle of $B(V_{1},V_{2})$. If $|N(u_1) \cap Y| =1$ and $|N(u_1) \cap X| =1$, we assume $\{x_1,y_1\} \subseteq N(u_1)$, then $xx_1u_1y_1yu_2x$ is a cycle of $B(V_{1},V_{2})$. If $|N(u_1) \cap X| =1$ and $|N(u_1) \cap Y| =0$, then $u_1u_2 \in E(G)$ and $xx_1u_1u_2x$ is a cycle of $B(V_{1},V_{2})$. Similarly, if $|N(u_1) \cap Y| =0$ and $|N(u_1) \cap X| =1$, then we can get a cycle of $B(V_{1},V_{2})$.
Therefore, $(V_1,V_2)$ is a 2-CIST-partition of $G$.  
 Thus, by Theorem \ref{thm1}, $G$ has two CISTs.

\vspace{0.5cm}

\noindent{\bf Subcase 1.2}  $d(u_2)=3$ and $u_1u_2 \notin E(G)$.

Let $V_{1}$ and $V_{2}$ be the same as in Subcase 1.1. By similar argument as in Subcase 1.1, we can get that $G[V_{1}]$, $G[V_{2}]$ and $B(V_{1},V_{2})$ are connected. For the existence of a cycle in $B(V_{1},V_{2})$, except for the last two cases in Subcase 1.1 which do not occur, everything else is the same. So $(V_1,V_2)$ is a 2-CIST-partition of $G$. Thus, by Theorem \ref{thm1}, $G$ has two CISTs.

\vspace{0.5cm}

\noindent{\bf Subcase 1.3}  $d(u_2)=3$ and $u_1u_2 \in E(G)$.

In this case, $N(u_2)=\{x,y,u_1\}$. Then $d(x_i,u_2)=2$ for any vertex $x_i \in X$ and $d(y_j,u_2)=2$ for any vertex $y_j \in Y$.  By Lemma \ref{lem1}, we have $N(x_i)=V(G) \setminus\{x_i,y,u_2\}$ for each vertex $x_i\in X$ and $N(y_j) =V(G) \setminus \{y_j,x,u_2\}$ for each vertex $y_j\in Y$. Thus, $N(u_1)=V(G)\setminus \{u_1\}$.

Let $V_{1}=\{y,u_1\}$, $V_{2}=V(G) \setminus V_1$. By the above discussion, $G[V_{1}]$, $G[V_{2}]$ and $B(V_{1},V_{2})$ are connected, and $yy_1u_1u_2y$ is a cycle of $B(V_{1},V_{2})$. Therefore, $(V_1,V_2)$ is a 2-CIST-partition of $G$. Thus, by Theorem \ref{thm1}, $G$ has two CISTs.

 \vspace{0.5cm}
\noindent{\bf Case 2}  $|M| \ge 3$.

\vspace{0.5cm}
The following claim is very useful to finish our proof of Theorem \ref{th2}.

\begin{cla}\label{cla3.4}
Each of the following holds:

\textnormal{(i)} $G$ is $m$-connected, where $m=|M| ;$ and

\textnormal{(ii)} if $U$ is an $m$-cut of $G$ containing both $x$ and $y$, then $G$ has two CISTs.
\end{cla}
	
\noindent
\textbf{Proof.} (i) Suppose not. Let $U$ be a minimum cut of $G$ with $|U| < m$ and let $G_1$, $G_2$ be two components of $G-U$. Since $U$ is minimum, there are two vertices $z_1 \in V(G_1)$ and $z_2 \in V(G_2)$ such that $d(z_1,z_2)=2$. By Lemma \ref{lem1}, $|N(z_1)\cap N(z_2)|\geq m$. This contradicts to $N(z_1)\cap N(z_2)\subseteq U$. Therefore, $G$ is $m$-connected.

(ii) Let $H_1$, $H_2$, $\dots$, $H_p$ ($p \ge 2$) be all the components of $G-U$. By (i), $U$ is a minimum cut of $G$. So there is at least a vertex $z_j \in V(H_j)$ such that $xz_j \in E(G)$ for each $j\in [p]$. By Lemma \ref{lem1}, $N(z_i) \cap N(z_j)=U$ for $1\leq i<j\leq p$. Since $\{x,y\} \subseteq U$, $\{z_1,z_2,\ldots,z_p\}\subseteq M$. Similarly, we can get that for each vertex $v\in V(G)\setminus U$, if $N_U(v)\neq \varnothing$, then $N_U(v)=U$ and $v\in M$. This implies that each vertex not in $M\cup U$ has no neighbor in $U$.

If there is a component $H_j$ such that $V(H_j)\setminus M\neq \varnothing$, then $N_{H_j}(U)\subseteq M\setminus \{z_1,z_2,\ldots,z_{j-1},$ $z_{j+1}, z_{j+2}, \ldots,z_{p}\}$ is a cut of $G$. This contradicts to (i). Therefore, $V(H_j)\subseteq M$ for each $j\in [p]$. Thus, $V(G)\setminus U\subseteq M$. So $n=|U|+|V(G)\setminus U|\le 2m$, $m \ge \lceil \frac{n}{2}\rceil \ge 4$.  

Without loss of generality, assume $\{u_1,u_2\}\subseteq V(G)\setminus U$. If $U \cap M = \varnothing$, then $V(G) \setminus U=M$. Choose any two vertices $a_1,a_2\in U \setminus \{x,y\}$. Let $W_1^1=\{x,y,u_1,u_2\}$, $W_2^1=V(G) \setminus W_1^1=(U\setminus \{x,y\})\cup (M\setminus \{u_1,u_2\})$. Obviously, $G[W_1^1]$ and $G[W_2^1]$ are connected. The induced subgraph of $B(W_1^1,W_2^1)$  by $\{x,y\}$ and $M \setminus \{u_1,u_2\}$ is connected and has a cycle $xu_3yu_4x$, and the induced subgraph of $B(W_1^1,W_2^1)$ by $\{u_1,u_2\}$ and $U \setminus \{x,y\}$ is also connected and has a cycle $u_1a_1u_2a_2u_1$. Therefore, $(W_1^1,W_2^1)$ is a 2-CIST-partition of $G$.  Thus, by Theorem \ref{thm1}, $G$ has two CISTs.

If $U \cap M \ne \varnothing$, then without loss of generality, assume $u_3\in U \cap M$. Let $W_1^2=\{x,u_1,u_3\}$, $W_2^2=V(G)\setminus W_1^2=(U\setminus \{x,u_3\})\cup (M\setminus \{u_1,u_3\})$. Obviously, $G[W_1^2]$, $G[W_2^2]$ and $B(W_1^2,W_2^2)$ are connected, and $xu_2u_3u_4x$ is a cycle of $B(W_1^2,W_2^2)$ if $u_4\notin U$, or $xu_2u_3yu_1u_4x$ is a cycle of $B(W_1^2,W_2^2)$ if $u_4\in U$. Therefore, $(W_1^2,W_2^2)$ is a 2-CIST-partition of $G$.  Thus, by Theorem \ref{thm1}, $G$ has two CISTs. $\hfill\Box$

\vspace{0.5cm} 
If  $D= \varnothing$, then let $V_1=\{x,y,u_1,\dots,u_{m-2}\}$ and $V_2=V(G) \setminus V_1$. Obviously, $G[V_1]$ is connected. If $V_1$ is an $m$-cut, then by Claim \ref{cla3.4} (ii), $G$ has two CISTs. If $V_1$ is not an $m$-cut, then $G[V_2]$ is connected. By Claim \ref{cla3.4} (i), $\delta(G)\geq m$. So every vertex of $\{u_1,\dots,u_{m-2}\}$ has at least one neighbor in $V_2$. Since the induced subgraph of $B(V_1,V_2)$ by $\{x,y\}$ and $V_2$ is connected, $B(V_1,V_2)$ is connected and has a cycle $xu_{m-1}yu_mx$. Therefore, $(V_1,V_2)$ is a 2-CIST-partition of $G$. Thus, by Theorem \ref{thm1}, $G$ has two CISTs.

Next, we assume that $D\neq \varnothing$.

\vspace{0.5cm} 
\noindent{\bf Subase 2.1}  $Y=\varnothing$.

Then $X=\varnothing$. Since $2|M|=d(x)+d(y)\geq n \ge 7$, $|M| \ge 4$. By Claim \ref{cla7}, we have $d(x,v_i)=2$ for any vertex $v_i \in D$. Thus, by the choice of $x$ and $y$, we have $N_M(v_i)=M$ and $D$ is an independent set of $G$. 

If $|D|\geq 2$, then let $V_1=\{x,u_1, u_2, v_1$\} and $V_2=V(G)\setminus V_1$. It is clear that $G[V_1]$ and  $G[V_2]$ are connected. Moreover, the induced subgraph of $B(V_1,V_2)$ by $\{x,v_1\}$ and $M\setminus \{u_1,u_2\}$ is connected and has a cycle $xu_3v_1u_4x$, and the induced subgraph of $B(V_1,V_2)$ by $\{u_1,u_2\}$ and $\{y\}\cup (D\setminus \{v_1\})$ is connected and has a cycle $yu_1v_2u_2y$. Therefore, $(V_1,V_2)$ is a 2-CIST-partition of $G$. Thus, by Theorem \ref{thm1}, $G$ has two CISTs.

If $|D|=1$, then we can claim that $M$ is not an independent set $G$, otherwise, $7\leq n\leq d(u_1)+d(u_2)\leq 6$, a contradiction. Suppose without loss of generality that $u_1u_2\in E(G)$. Let $V_1=\{x,u_1,v_1\}$ and $V_2=V(G)\setminus V_1$. It is clear that $G[V_1]$, $G[V_2]$ and $B(V_1,V_2)$ are connected, and $xu_2v_1u_3x$ is a cycle of $B(V_1,V_2)$. Therefore, $(V_1,V_2)$ is a 2-CIST-partition of $G$. Thus, by Theorem \ref{thm1}, $G$ has two CISTs.

\vspace{0.5cm} 
\noindent{\bf Subase 2.2}  $Y\neq \varnothing$.
\vspace{0.5cm}

By Claim \ref{cla7}, for any vertex $v_i \in D$, we have $d(y,v_i)=2$ and $|N(y) \cap N(v_i)| \ge m$.  
Select a subset $S \subseteq M \cup Y$ such that 

(C3) $N_S(v_i) \ne \varnothing$ for any vertex $v_i \in D$,  subject to (C2);

(C4) $|S|$ is as small as possible, subject to (C3);

(C5) $|S \cap M|$ is as large as possible, subject to (C4).

By the choice of $S$, we can get the following fact:

\begin{fa} \label{fa3.1}
Each of the following holds:

$(1)$ for any two vertices $u,v\in S$, we have $N_D(u)\setminus N_D(v)\neq \varnothing$ and $N_D(v)\setminus N_D(u)\neq \varnothing;$

$(2)$ for any subset $S'\subseteq S$, $|S'|\leq |N_D(S')|.$  Particularly, $|S|\leq |D| \le m-2$.
 
\end{fa}

\noindent
\textbf{Proof.} (1) Suppose not. There are two vertices $u,v \in S$ such that $N_D(u)\subseteq N_D(v)$. Then $S'=S\setminus \{u\}$ is a subset of $M\cup Y$ satisfying (C3) and $|S'|<|S|$, which contradicts to (C4).

(2) Let $S'=\{s_1,\dots,s_l\}\subseteq S$. By (1), there exist distinct vertices $v_1,v_2,\dots, v_l\in D$ such that $s_iv_i \in E(G)$ for $i \in [l]$. Let $L=\{s_1v_1,s_2v_2,\ldots, s_lv_l\}$. Then, $|S'|=|N_L(S')|\leq |N_D(S')|$.
Particularly,  $|S| \le |N_D(S)|=|D| \le m-2$. $\hfill\Box$

\vspace{0.5cm}	
\noindent{\bf Subcase 2.2.1}\label{sub2.2.1} $S \cap M \ne \varnothing$. 

In this case, let $V_{1}^{1}=\{x,y\} \cup S$, $V_{2}^{1}=V(G) \setminus V_1^1$. Since $S \cap M \ne \varnothing$, $G[V_1^1]$ is connected. If $G[V_2^1]$ is not connected, then by Claim \ref{cla3.4} (i), $V_1^1$ is an $m$-cut of $G$ containing both $x$ and $y$. Hence, by Claim \ref{cla3.4} (ii), $G$ has two CISTs. 

Next, we assume that $G[V_2^1]$ is connected. We will prove that there is no tree component in $B(V_{1}^{1},V_{2}^{1})$. Therefore, $(V_1^1,V_2^1)$ is a 2-CIST-partition of $G$. Thus, by Theorem \ref{thm1}, $G$ has two CISTs.

\begin{prop}
There is no tree component in $B(V_{1}^{1},V_{2}^{1})$.

\end{prop}
	
\noindent
\textbf{Proof.} By Fact \ref{fa3.1},  $|M\setminus S| \ge 2$. Without loss of generality, assume that $\{u_1,u_2\} \subseteq M\setminus S$. Clearly, the induced subgraph of $B(V_{1}^{1},V_{2}^{1})$ by $\{x,y\}$ and $(X \cup M \cup Y)\setminus S$, say $B_3$, is connected and has a cycle $xu_1yu_2x$. 

Next, we will prove that $B(V_{1}^{1},V_{2}^{1})$ is connected. Otherwise, let $G_1$, $G_2$, $\dots$, $G_q$ ($q \ge 2$)  be all the components of $B(V_{1}^{1},V_{2}^{1})$. Without loss of generality, we assume $B_3\subseteq G_1$. 
Let $S_i=V(G_i) \cap S$ and $D_i=V(G_i) \cap D$ for $i\in [q]$. By Fact \ref{fa3.1}(2), $|S_i|\leq |D_i|$ for each $i\in [q]$. Hence, $|D|-|S|=\sum\limits_{i=1}^q (|D_i|-|S_i|)\geq |D_i|-|S_i|$, that is $|D_i|\leq |D|-|S|+|S_i|$ for $i\in [q]$.  

We can claim that $|S_i|=1$ for $2\leq i\leq q$. Otherwise, without loss of generality, suppose $|S_2| \ge 2$ and $s_1,s_2\in S_2$. For convenience, let $D'=N_D(s_1) \cup N_D(s_2)$. 
By Fact \ref{fa3.1} (2), we have $|D \setminus D'|\geq |S\setminus \{s_1,s_2\}|=|S|-2$. Hence, $|D'|=|D|-|D \setminus D'| \le |D|-|S|+2$. By Fact \ref{fa3.1} (1), there is a vertex $v'\in N_D(s_1)\setminus N_D(s_2)$ and a vertex $v''\in N_D(s_2)\setminus N_D(s_1)$. Since $ys_1\in E(G)$ and $yu_1 \in E(G)$, $d(u_1,s_1)=2$. Thus $|N(u_1) \cap N(s_1)|\ge m$. Since $N(u_1) \cap N(s_1) \subseteq (V_1^1 \setminus S_2) \cup (D' \setminus \{v''\})$, $|N(v') \cap N(s_1)| \le |V_1^1 \setminus S_2|+|D' \setminus \{v''\}|\leq (|S|+2-|S_2|)+(|D|-|S|+1)\leq |D|+1<m$, a contradiction. Therefore, $|S_i|=1$ for $2 \le i \le q$. Let $S_i=\{s_i\}$ and $v_i\in D_i$ for  $2 \le i \le q$.

If $q\geq 3$, then $N(u_1) \cap N(s_2) \subseteq (V_1^1 \setminus \{s_2,s_3\}) \cup D_2$. Thus, $m\leq |N(u_1) \cap N(s_2)|\leq |V_1^1 \setminus \{s_2,s_3\}|+|D_2| \leq |S|+(|D|-|S|+1)=|D|+1<m$, a contradiction. Thus $q=2$.

Since $d(u_1,s_2)=2$ and $N(u_1) \cap N(s_2) \subseteq (V_1^1 \setminus \{s_2\}) \cup D_2$, we have $m\leq |N(u_1)\cap N(s_2)|\leq (|S|+1)+|D_2|\leq (|S|+1)+(|D|-|S|+1)\leq |D|+2\leq m$. Therefore, $N(u_1)\cap N(s_2)=(V_1^1\setminus \{s_2\})\cup D_2$.
Then $s_2\in M$. 

If $|S|\geq 2$, then let $s_1\in S_1$ and $v_1\in D_1$. Since $s_2s_1\in E(G)$, $s_2v_1\notin E(G)$ and $N(v_1)\cap N(s_2)\subseteq (V_1^1\setminus \{x,y,s_2\})\cup D_2$, we have $m\leq |N(v_1)\cap N(s_2)|\leq (|S|-1)+|D_2|\leq (|S|-1)+(|D|-|S|+1)\leq |D|<m$, a contradiction.

If $|S|=1$, then $d(s_2,y_1)=2$ and $|N(y_1)\cap N(s_2)|\geq m$. However, $N(y_1) \cap N(s_2) \subseteq (V_1^1\setminus \{s_2,x\}) \cup D_2$, and so $|N(y_1)\cap N(a)|\leq |S|+(|D|-|S|+1)\leq |D|+1< m$, a contradiction. 

\begin{figure}[htbp]  
  \centering 
  \includegraphics[width=0.4\textwidth]{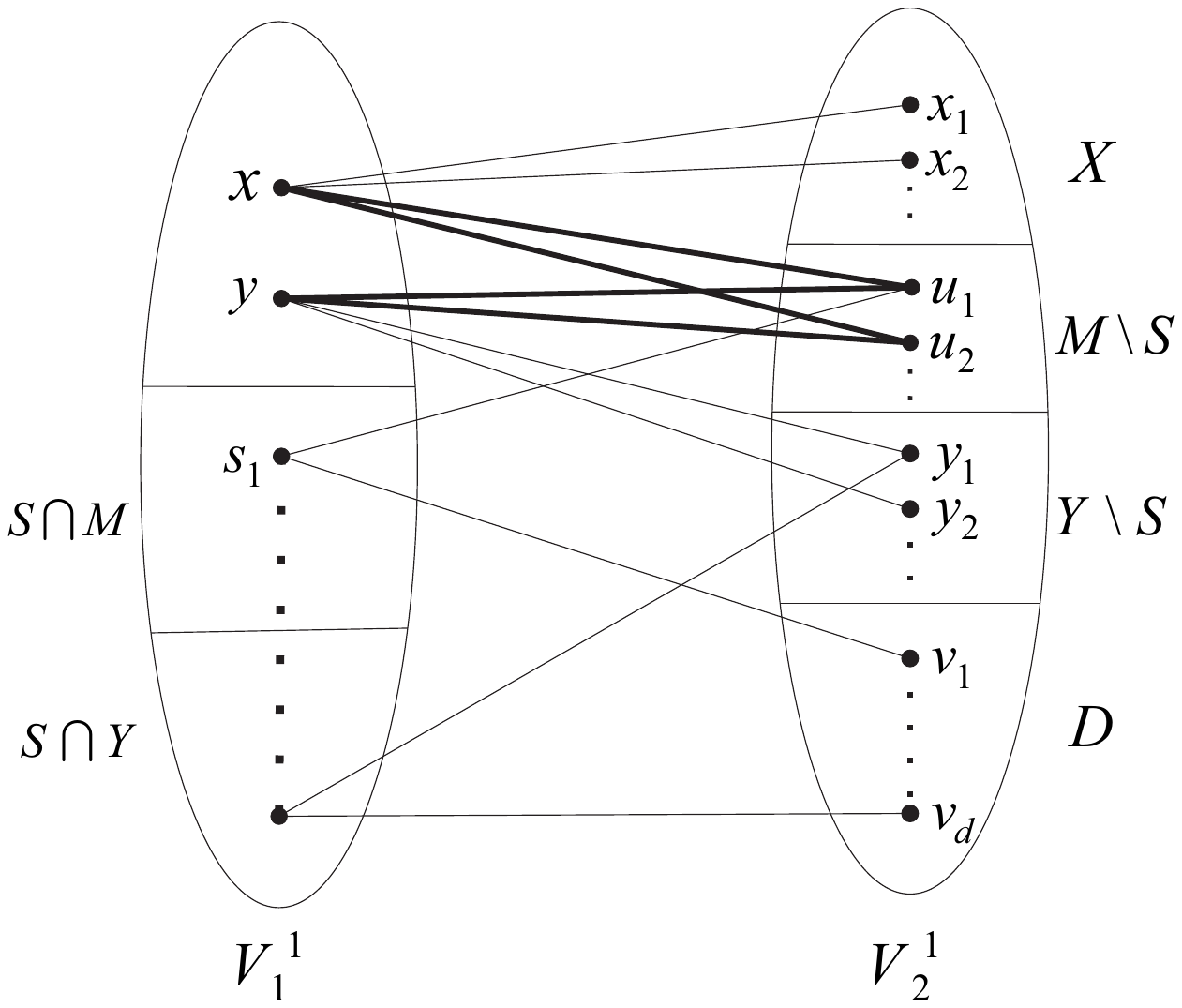} 
  \vspace{-0.7cm}
  \caption{\footnotesize The illustration of a 2-CIST-partition in Subcase 2.2.1.}  
  \label{fig:fig2}  
\end{figure}

Therefore, $B(V_{1}^{1},V_{2}^{1})$ is connected and has no tree component (see Figure \ref{fig:fig2}). $\hfill\Box$

\vspace{0.5cm}

\noindent{\bf Subcase 2.2.2} $S \cap M = \varnothing$.

In this case, $S\subseteq Y$. 
Set $R_j=\{v\in D|u_jv\notin E(G), u_j\in M\}$ for $j\in [m]$. We can easily get that $|R_j|\geq |S|$ for each $j\in [m]$. Otherwise, by Fact \ref{fa3.1} (1), there is a vertex $y_i\in S$ and a vertex $u_j\in M$ such that $N_D(y_i)\subseteq N_D(u_j)$. Then $S'=(S \setminus \{y_i\}) \cup \{u_j\} \subseteq M \cup Y$ is a set satisfying (C3) and (C4), but $|S' \cap M|>|S \cap M|$, which contradicts to (C5). 

\begin{cla}\label{cla3.6}
Each vertex of $S$ has at least $|S|+1$ neighbors in $X\cup Y\cup M$.
\end{cla}
	
\noindent
\textbf{Proof.} Suppose to the contrary that there is a vertex $u\in S$ such that $|N_{X\cup Y\cup M}(u)|\leq |S|$. Since $|M|>|S|$ and $N(y)=M\cup Y$, there is a vertex $u_j\in M$ such that $d(u,u_j)=2$. Then $|N(u)\cap N(u_j)|\geq m$. On the other hand, since $N(u)\cap N(u_j)\subseteq N_{X\cup Y\cup M}(u)\cup \{y\}\cup N_D(u_j)$, we can get that $|N(u)\cap N(u_j)|\leq |S|+1+(|D|-|R_j|)\leq |S|+1+(|D|-|S|)=|D|+1<m$, a contradiction. $\hfill\Box$

\vspace{0.5cm}
Let $V_{1}^2=\{x,y,u_1\} \cup S$ and $V_{2}^2=V(G) \setminus V_1^2$. Clearly, $G[V_{1}]$ is connected, the induced subgraph of $B(V_{1}^2,V_{2}^2)$ by $\{x,y\}$ and $X\cup (Y\setminus S)\cup (M\setminus \{u_1\})$ is connected and has a cycle $xu_2yu_3x$. By Claim \ref{cla3.6}, each vertex of $S$ has a neighbor in $X\cup (Y\setminus S)\cup (M\setminus \{u_1\})$. Therefore, the subgraph of $B(V_{1}^2,V_{2}^2)$ induced by $S\cup \{x,y\}$ and $V_2^2$ is connected.

If $|S|\leq m-3$, then $|V_{1}^2|\leq m$. If  $V_{1}^2$ is an $m$-cut of $G$, then by Claim \ref{cla3.4} (ii), $G$ has two CISTs. Otherwise, by Fact \ref{fa3.1} (1), $G[V_{2}^2]$ is also connected. Since $d(u_1)\geq m>|V_1^1\setminus \{u_1\}|$, $u_1$ has a neighbor   in $V_2^2$. Therefore, $B(V_{1}^2,V_{2}^2)$ is connected. This implies that $(V_1^2,V_2^2)$ is a 2-CIST-partition of $G$. Thus, by Theorem \ref{thm1}, $G$ has two CISTs.

Next, we assume $|S|=m-2$, then $|V_{1}^2|=m+1$. By Fact \ref{fa3.1}(2), $d\le m-2=|S|\le |D|=d$. So $|S|=|D|=d=m-2$. By the choice of $S$, $B(S,D)$ is isomorphism to $(m-2)K_2$ and $N_M(v_i)=\varnothing$ for every vertex $v_i\in D$. Without loss of generality, assume $S=\{y_1, \dots, y_{m-2}\}$ and $y_iv_i\in E(G)$, $i \in [m-2]$.

If there is a vertex $y_i\in S$ such that $V_{1}^2\setminus \{y_i\}$ is an $m$-cut of $G$, then by Claim \ref{cla3.4} (ii), $G$ has two CISTs. Suppose $V_{1}^2\setminus \{y_i\}$ is not an $m$-cut of $G$ for each vertex $y_i\in S$. We will prove that $(V_1^2,V_2^2)$ is also a 2-CIST-partition of $G$. 

\begin{prop}\label{pro3.2}
$G[V_{2}^{2}]$ is connected.
\end{prop}
	
\noindent
\textbf{Proof.} Suppose not. Let $G_1$, $G_2$, $\dots$, $G_q$ ($q \ge 2$) be all the components of $G[V_2^2]=G-V_1^2$.  
Without loss of generality, suppose $v_1\in V(G_1)\cap D$. 
Since $V_1^2\setminus \{y_1\}$ is not a cut, but $V_1^2$ is a cut of $G$, we can get that $y_1$ is a cut-vertex of $G-(V_1^2\setminus \{y_1\})$. Then there is at least one vertex $z_i\in V(G_i)$ ($i \in [q]$) such that $y_1z_i\in E(G)$. Since $v_1y_1\in E(G)$, $d(v_1,z_2)=2$. So $|N(v_1)\cap N(z_2)|\geq m\geq 3$. However, $N(v_1)\cap N(z_2)= \{y_1\}$, a contradiction. Therefore, $G[V_{2}^{1}]$ is connected. $\hfill\Box$

\begin{prop}\label{pro3.3}
There is no tree component in $B(V_{1}^{2},V_{2}^{2})$.
\end{prop}

\begin{figure}[htbp]  
  \centering 
  \includegraphics[width=0.4\textwidth]{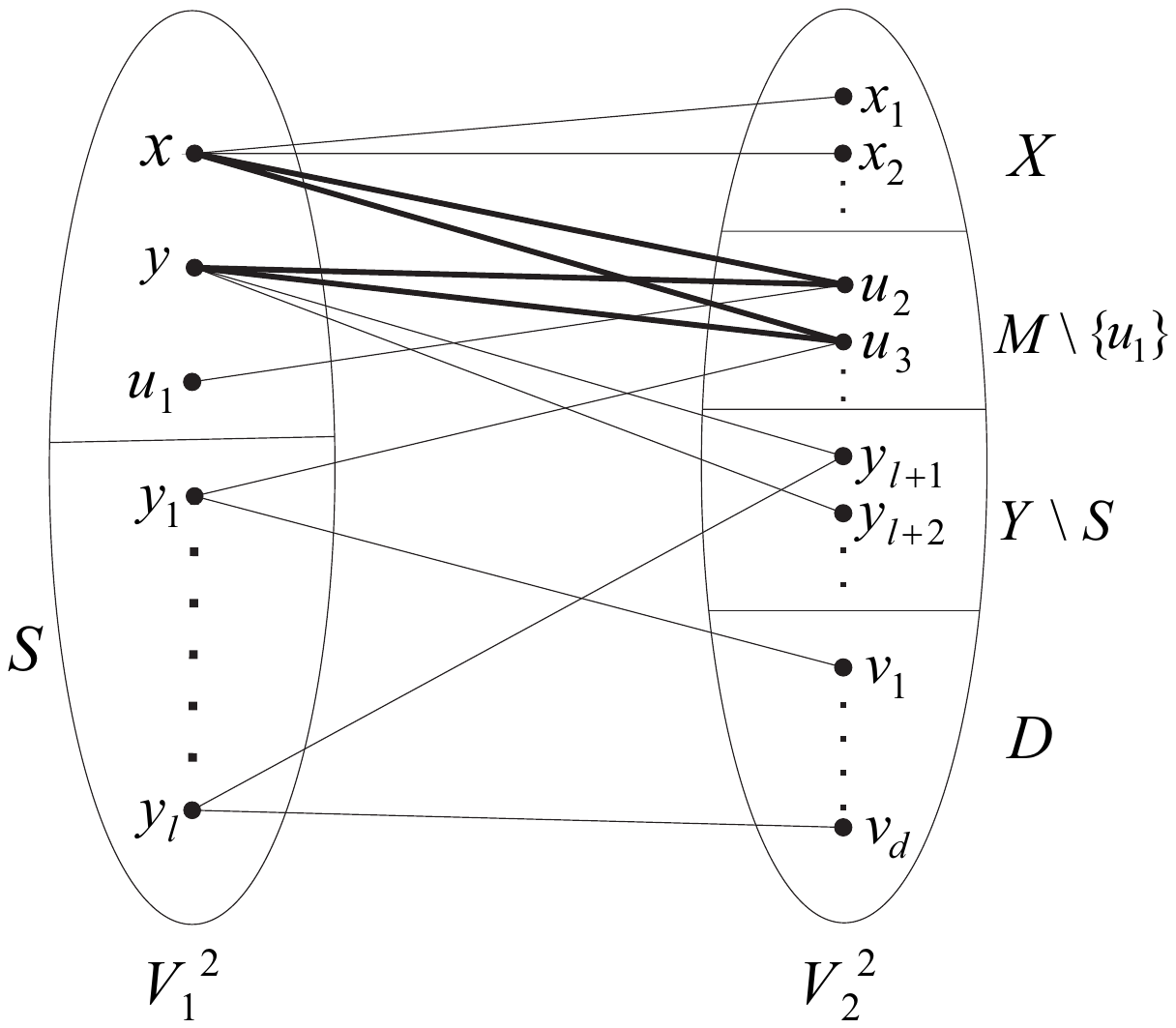} 
  \vspace{-0.7cm}
  \caption{\footnotesize The illustration of a 2-CIST-partition in Subcase 2.2.2.}  
  \label{fig:3}  
\end{figure}

\noindent
\textbf{Proof.} Since $N_M(v_i)=\varnothing$ for every vertex $v_i\in D$ and $d(v_i)\geq m$, $|Y|\geq m$. So $Y\setminus S\neq \varnothing$. If $u_1$ has no neighbor in $V_2^2$, then $d(u_1, y_{m-1})=2$ and thus $|N(u_1)\cap N(y_{m-1})|\geq m$. On the other hand, since $N(u_1)\cap N(y_{m-1})\subseteq V_1^2\setminus \{x,u_1\}$, $|N(u_1)\cap N(y_{m-1})|\leq |V_1^2|-2=m-1$, a contradiction. Therefore, $B(V_{1}^{2},V_{2}^{2})$ is connected and has a cycle $xu_2yu_3x$ (see Figure \ref{fig:3}). Therefore, there is no tree component in $B(V_{1}^{2},V_{2}^{2})$. $\hfill\Box$

Therefore, by Propositions \ref{pro3.2} and \ref{pro3.3}, $(V_1^2,V_2^2)$ is a 2-CIST-partition of $G$ (see Figure 3).  Thus, by Theorem \ref{thm1}, $G$ has two CISTs.

Now, we have completed the proof of Theorem \ref{th2}. $\hfill\Box$

\section*{Declaration of competing interest}
The authors declare that they have no known competing financial interests or personal relationship that could have appeared to influence the work reported in this paper.

\section*{Acknowledgment}
This work was supported by National Natural Science Foundation of China (No.12371356).

\end{document}